\documentclass[a4paper,10pt]{article}
\usepackage[cp1251]{inputenc}
\usepackage[english]{babel}

\usepackage{amssymb}\usepackage{amsmath}
\usepackage{epsfig}
\usepackage{graphics}
\def\bl{\rule[-1mm]{2.4mm}{2.4mm}}

\def\be{\begin{equation}}
\def\ee{\end{equation}}

\newtheorem{thrm}{\bf Theorem}
\newtheorem{lmm}{\bf Lemma}

\newtheorem{rmk}{\bf Remark}

\begin{document}

\title {Capacity of several aligned segments}
\author{\copyright 2015 ~~~~A.B.Bogatyrev, ~~O.A.Grigoriev
\thanks{Both autors are supported by RFBR grant 13-01-00115, first author is also supported by 13-01-12417ofi-m2 and RAS Program
"Modern problems of theoretical mathematics"}}
\date{}
\maketitle

\begin{abstract}\bf
In this note we present a universal formula in terms of theta functions for the Log- capacity of several segments on a line.
The case of two segments was studied by N.I.Akhiezer (1930);  three segments were considered by A.Sebbar and T.Falliero (2001).
\end{abstract}

For physicist the capacity is a coefficient relating the charge and the voltage of a condenser, or equivalently  the energy of an electric field inside the condenser under the unit voltage. For mathematician the capacity is a certain property of a set which shows how 'massive' is it. This concept turned out to be very usefull in approximation theory, geometric function theory \cite{Go}, partial differential equations, potential theory to name a few. For a compact set in the complex plane its capacity coinsides with the Chebyshev constant, the transfinite diameter, conformal radius (for the simply connected sets) and a simple formula links it to the Robin constant of the set. 

Let $E$ be a collection of $g+1$ segments on a real line:
\be
\label{Edef}
E:=\cup_{j=0}^g ~[e_{2j+1},e_{2j+2}]
\ee
with strictly increasing sequence of endpoints $e_j$. We know that the capacity of a set 
is independent of its  translations and homogenious with respect to dilations, so 
without loss of generality we set the most left point of $E$ being $0$ and the most right point being $1$:
$$
e_1=0; \qquad e_{2g+2}=1.
$$ 
The capacity $C:=Cap(E)$  is defined in terms of the asymptotic of the Green's function $G_E(x)$ of this set at infinity:
\be
\label{CapE}
G_E(x)=\log|x| -\log C +~o(1), \qquad x\to \infty.
\ee
We are going to deduce a closed formula for $Cap(E)$ in terms of Riemann's theta functions so 
we have to introduce some necessary (but standard) constructions related to Riemann surfaces.

\section{Riemann surface and theta functions}
The Riemann surface associated to our problem is the double cover of the sphere ramified over the endpoints of the segments from $E$. This is a genus $g$ compact surface $\cal X$ with its affine part given by the equation  
\be
\label{X}
w^2=\prod_{j=1}^{2g+2}(x-e_j).
\ee
The surface admits the hyperelliptic involution $J(x,w):=(x,-w)$
with $2g+2$ fixed points $P_s=(e_s,0)$ as well as an anticonformal involution (reflection)
$\bar{J}(x,w):=(\bar{x},\bar{w})$. Fixed points of the reflection $\bar{J}$ make up the  \emph{real ovals} of the surface;
fixed points of another anticonformal involution $J\bar{J}$ are known as coreal ovals. Coreal ovals cover the set $E$; real ovals cover the compliment
$\hat{\mathbb R}\setminus E$, in particular two points $\infty_\pm:=(+\infty,\pm\infty)$
from the same real oval cover the infinity.

\begin{figure}[t]
\begin{picture}(150,55)
\put(0,0){\includegraphics[scale =.75]{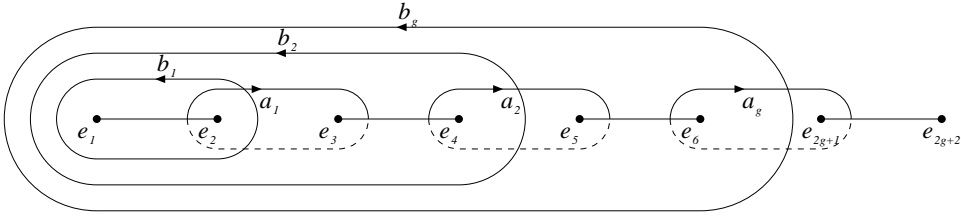}}
\end{picture}
\caption{\small Symplectic basis in the 1-homologies of the curve $\cal X$}
\label{Basis}
\end{figure}

We introduce a symplectic basis in the homologies of $\cal X$ as shown in the Fig. \ref{Basis}.
Dual basis of holomorphic differentials satisfies the normalization conditions
\be
\label{duNorm}
\int_{a_j}du_s:=\delta_{js},
\ee
and generates the period matrix
\be
\label{Pi}
\int_{b_j}du_s=:\Pi_{js}.
\ee

One easily checks that the introduced bases of cycles and differentials behave as follows under the reflection:
\be
\bar{J}a_s=a_s;
\qquad  \bar{J}b_s=-b_s;
\ee
\be
\label{du} 
\qquad \bar{J}du_s=\overline{du_s}; 
\qquad s=1,\dots,g.
\ee
Cycles surviving after the reflection are called \emph{even} and those changing sign are called \emph{odd}. Differentials with the property (\ref{du}) are usually called \emph{real} and their periods along even/odd cycles are real/pure imaginary respectively, see \cite{B1, B2} for more details. In paricular, our period matrix $\Pi$ is purely imaginary.
Standard facts are the symmetry of this matrix and positive definiteness of the imaginary part \cite{FK}.

Given the period matrix, we define the complex $g$-dimensinal torus known as a 
Jacobian of the complex curve $\cal X$:  
\be
Jac({\cal X}):=\mathbb{C}^g/L(\Pi), 
\qquad\qquad L(\Pi):=\mathbb{Z}^g+\Pi\mathbb{Z}^g.
\label{Jac}
\ee
The Abel-Jacobi (AJ) map embeds the curve to its Jacobian
\be
\label{AJmap}
u(P):=\int_{P_1}^P du \bmod L(\Pi) \in Jac({\cal X}),
\qquad du:=(du_1,du_2,\dots,du_g)^t.
\ee 

It is convenient to represent points  $u\in\mathbb{C}^g$ as the so called theta characteristics, i.e. couple of real g-vector columns $\epsilon, \epsilon'$:
\be
u=(\epsilon'+\Pi\epsilon)/2.
\ee
The points of Jacobian in this notation correspond to two vectors with real entries modulo 2. Second order points of Jacobian are
$2\times g$ matrices with binary entries. In particular, the images of the hyperelliptic involution fixed points $P_s$ the  under the AJ map (\ref{AJmap}) are as follows:

\begin{tabular}{c|c|c}
$P_s$ & $u(P_s)~ mod~ L(\Pi)$ & $[\epsilon, \epsilon']^t$\\
\hline
$P_1$&$0$&$\tiny
\left[\begin{array}{c}00\dots0\\00\dots0\end{array}\right]$\\
$P_2$ & $\Pi_1/2$ & $\tiny
\left[\begin{array}{c} 100\dots\\000\dots\end{array}\right]$\\
$P_3$&$(\Pi_1+E_1)/2$&$\tiny
\left[\begin{array}{c} 100\dots\\100\dots\end{array}\right]$\\
$P_4$&$(\Pi_2+E_1)/2$&$\tiny
\left[\begin{array}{c} 010\dots\\100\dots\end{array}\right]$\\
$P_5$&$(\Pi_2+E_1+E_2)/2$&$\tiny
\left[\begin{array}{c}010\dots\\110\dots\end{array}\right]$\\
$P_6$&$(\Pi_3+E_1+E_2)/2$&$\tiny
\left[\begin{array}{c}0010\dots\\1100\dots\end{array}\right]$\\
\vdots&\vdots&\vdots\\
$P_{2g}$&$(\Pi_g+E_1+E_2+\dots+E_{g-1})/2$&$\tiny
\left[\begin{array}{c} 0\dots01\\1\dots10\end{array}\right]$\\
$P_{2g+1}$&$(\Pi_g+E_1+E_2+\dots+E_g)/2$&$\tiny
\left[\begin{array}{c} 0\dots01\\11\dots1\end{array}\right]$\\
$P_{2g+2}$&$(E_1+E_2+\dots+E_g)/2$&$\tiny
\left[\begin{array}{c}00\dots0\\11\dots1\end{array}\right]$\\
\label{AJPj}
\end{tabular}   

Here $E_s$ and $\Pi_s$ are the columns of the identity matrix and the period matrix respectively.
Please note that we use a nonstandard notation of theta characteristic as two column vectors written one after another. Usually the transposed matrix is used.

The following absolutely convergent series is known as theta function:
\be
\label{thetadef}
\theta(u, \Pi):=\sum\limits_{m\in\mathbb{Z}^g}
\exp(2\pi i m^tu+\pi i m^t\Pi m),
\qquad u\in\mathbb{C}^g;\Pi=\Pi^t\in\mathbb{C}^{g\times g}; Im~\Pi>0.
\ee
Often it is convenient to consider theta function with charachteristics, a slight modification of the latter:

\be
\label{thetachardef}
\theta[2\epsilon, 2\epsilon'](u, \Pi):=\sum\limits_{m\in\mathbb{Z}^g}
\exp(2\pi i(m+\epsilon)^t(u+\epsilon')+\pi i (m+\epsilon)^t\Pi (m+\epsilon))
\ee

$$
=\exp(i\pi\epsilon^t\Pi\epsilon+2i\pi\epsilon^t(u+\epsilon'))
\theta(u+\Pi\epsilon+\epsilon',\Pi),
\qquad \epsilon,\epsilon'\in\mathbb{R}^g.
$$
Matrix argument $\Pi$ of theta function may be omitted if it does not lead to a confusion.
Vector argument $u$ may also be ommited which mean that we set $u=0$, the value of the theta function is 
called the theta constant in this case (note, there remains the dependence on the period matrix) 

This function (\ref{thetadef}) has the following easily checked quasi-periodicity properties with respect to the lattice $L(\Pi)$:
\be
\label{quasiperiod}
\theta(u+m'+\Pi m; \Pi)=\exp(-i\pi m^t\Pi m-2i\pi m^tu)\theta(u;\Pi),
\qquad m,m'\in\mathbb{Z}^g.
\ee
Theta functions with charachteristics have similar transformation rules \cite{Mu,RF} which may be easily derived from the above formula
and we omit them for brevity.

\begin{rmk}
\label{RemTheta}
(i) Theta function with integer characteristics $[\epsilon, \epsilon']$ is either even or odd depending on the 
parity of the inner product $\epsilon^t\cdot\epsilon'$. In particular, all odd theta constants are zeroes.

(ii) Adding matrix with even entries to integer theta characteristics can at most spoil the sign of the theta function. 
Hence, the binary arithmetic plays a great role in the calculus of theta functions.
\end{rmk}

Theta function may be considered as a multivalued function in the Jacobian
or as a section of a certain line bundle. The zero set of theta function -- the theta divisor -- is well defined 
in the Jacobian since the  exponential factor in the right hand side of (\ref{quasiperiod}) does not vanish.
The theta divisor is described by so called Riemann vanishing theorems  \cite{RF,FK}.
One of the important ingredients in those theorem is a vector of Riemann's constants $\cal K$
which depends on the choice of homology basis and the initial point in AJ map. 
In the above setting the vector of Riemann's constants may be found by a straightforward computation \cite{Mu} or by 
some combinatorial argument \cite{FK} and corresponds to a characteristic 
$$
{\cal K} \sim \left[\begin{array}{c}\dots11111\\ \dots10101\end{array}\right]^t.
$$

Zeroes of theta function transferred to the surface by the AJ map are described by the following 

\begin{thrm}[Riemann] Let $D_g$ be a degree $g$ positive non-special divisor on the curve $\cal X$, then
the function $\theta(u(P)-{\cal K}-u(D_g))$ of the argument $P\in{\cal X}$
has exactly $g$ zeroes at the points of $D_g$ counting their multiplicities.
\label{ThR}
\end{thrm}

Here and in what follows we use the standard notions of the function theory on Riemann surfaces \cite{FK,RF,Mu}.
\emph{Divisor} is a finite set of points on a surface taken with integer weights: $D=\sum_j m_jP_j$, $m_j\in\mathbb{Z}$, $P_j\in{\cal X}$. The \emph{degree of a divisor} is the sum of all weights of participating points: $\deg D:=\sum_j m_j$. 
\emph{Positive divisor} is the divisor with positive weights only.
The notion of divisor is very natural and useful in e.g. description of zeroes/poles of  meromorphic functions or differentials on a surface.

{\emph Index of speciality} $i(D)$ of a positive divisor $D$ is the dimension of the space of holomorphic differentials
vanishing at the points of $D$ (counting multiplicities). Divisors with index $i(D)>g-\deg D$ are called {\emph special}, which means that their points are not in generic positions.  Two divisors are {\emph linearly equivalent} iff their
difference is a set of zeroes and poles (latter should be taken with negative multiplicities) of a meromorphic function.
Riemann-Roch theorem \cite{FK} says that class of a positive divisor contains a unique element iff this divisor is not special.

\section{Green's function}
The goal of this section is to give a closed expression for the Green's function for the collection $E$ of aligned segments.
There exists a unique 3rd kind abelian differential $d\eta$ with simple poles at two distinguished points $\infty_\pm$ of the curve  $\cal X$ with residues $\mp1$ respectively and purely imaginary periods. Normalization conditions for this differential imply that it is real and therefore it's $a-$ periods (along even cycles) vanish \cite{B2}.

\begin{lmm} Green's function for $E$ has the representation:
\be
\label{Green}
G_E(x)=|Re(\int_Q^P d\eta)|,
\ee
where the upper limit $P$ covers the argument $x$ of Green's function, 
$x(P)=x$ and the lower limit $Q$ lies on the coreal oval, $x(Q)\in E$.

\end{lmm}
Proof.
We lift the compliment $\hat{\mathbb{C}}\setminus E$ to the Riemann surface (\ref{X})
so that infinity is mapped to the distinguished point $\infty_+$ of $\cal X$, and call it the top sheet of the surface. 
Other lifting we call the lower sheet. The function $G_E(x)$ on the top sheet of the surface and $-G_E(x)$ on the lower sheet make up a harmonic function with two log poles at $\infty_\pm$. It is a real part of some 3rd kind abelian integral which has the same normalizations as the introduced above integral $\eta=\int d\eta$. It remains to check that integral of $d\eta$ along any coreal oval is purely imaginary. This holds since the distinguished differential is real. ~~\bl

The explicit formula for the 3rd kind  abelian integrals in terms of theta functions was given yet by Riemann \cite{FK,RF}:
\begin{thrm}
$$
G_E(x)=\left|
\log
\left|
\frac{\theta[\epsilon,\epsilon'](u(P)+u(\infty_+))}{\theta[\epsilon,\epsilon'](u(P)-u(\infty_+))}
\right|
\right|,
$$
where $x(P)=x$ and $[\epsilon,\epsilon']$ is the integer (odd non-singular) theta charachteristic corresponding to half-period ${\cal K}+u(D_{g-1})$ where  $D_{g-1}$ is a sum of any $g-1$ different branch points:
$D_{g-1}=\sum_{s\in I}P_s$, $I\subset\{1,2,\dots, 2g+2\}$; $\#I=g-1$.
\end{thrm} 
Proof is based in Riemann's theorem on the zeros of theta functions:
the divisors $D_g^\pm:=D_{g-1}+\infty_\pm$ are non-special, therefore the fraction 
$\theta(u(P)-{\cal K}-u(D_g^-)/\theta(u(P)-{\cal K}-u(D_g^+)$ as a multivalued function of $P\in\cal X$
has a unique zero at $P=\infty_-$ and a unique pole at $P=\infty_+$. Same zero and pole has the function $\exp(\int^P d\eta)$, moreover both functions acquire the same factors when the argument $P$ goes around the cycles of the surface. 
Latter claim follows from the transformation rules of theta fuctions (\ref{quasiperiod}) and Riemann bilinear relations \cite{FK}.
From Liouville's theorem in now follows that both functions differ by a constant nonzero factor. The usage of theta charachteristics allows us to simplify the expression for the abelian integral. ~~\bl

\section{Independent variable and the Main formula}
Formula (\ref{CapE}) for the capacity  contains yet another ingredient along with the Green's function: the independent variable $x$. The representation for the hyperelliptic projection of the surface (\ref{X})  to the sphere in terms of variables of Jacobian  is well known \cite{FK}:

\begin{thrm}
$$
x(P)=\frac{\theta^2[\epsilon,\epsilon'](u(P))}{\prod_\pm\theta[\epsilon,\epsilon'](u(P)\pm u(\infty_+))}
\frac{\theta^2[\xi,\xi'](u(\infty_+))}{\theta^2[\xi,\xi']}
$$
where $[\epsilon,\epsilon']$ is again the integer theta charachteristic corresponding to half-period ${\cal K}+u(D_{g-1})$ and
$D_{g-1}$ is a positive divisor of any $g-1$ different branch points with the exception of first and last ones;
$[\xi,\xi']=[\epsilon,\epsilon']+[(0\dots0)^t,(1\dots1)^t] \bmod 2$
\end{thrm}
Proof. The function $x(P)$ on the surface $\cal X$ has double pole at $P=P_1$ and two poles at $P=\infty_\pm$, same as the fraction 
$\theta^2(u(P)-{\cal K}-u(D_{g-1}))/\prod_\pm\theta(u(P)-{\cal K}-u(D_{g-1})\pm u(\infty_+))$  -- it follows from Riemann's theorem on zeroes of thetas. The latter expression is single -valued on the surface and therefore differs by a constant factor from  the projection $x(P)$. This factor may be found if we evaluate both functions at $P=P_{2g+2}$. The resulting expression may be simplified by the usage of theta charachteristics.  ~~\bl

Let us assemble the previously obtained formulas:
$$
Cap(E)=\lim_{x\to\infty}|x|\exp(-G_E(x))=\lim_{P\to\infty_\pm}|x(P)|
\left|\frac{\theta[\epsilon,\epsilon'](u(P)\mp u(\infty_+))}{\theta[\epsilon,\epsilon'](u(P)\pm u(\infty_+))}\right|=
$$
\be
\label{main}
=\left|
\frac{\theta[\epsilon,\epsilon'](u(\infty_+)\theta[\xi,\xi'](u(\infty_+)}
{\theta[\epsilon,\epsilon'](2u(\infty_+))\theta[\xi,\xi'](0)}
\right|^2.
\ee
where $[\epsilon,\epsilon']$ is the integer theta charachteristic corresponding to half-period ${\cal K}+u(D_{g-1})$ with
$D_{g-1}$ a positive divisor of any $g-1$ different branch points, except of the first and the last one;
$[\xi,\xi']=[\epsilon,\epsilon']+[(0\dots0)^t,(1\dots1)^t] \mod 2$.

To make this formula computationally effective one has to calculate the period matrix $\Pi$ and the image of infinity in the Jacobian $u(\infty_+)$. However, there is a way to avoid usage of quadrature rules -- see e.g. techniques developed in \cite{B,G,B3}. The latter approach is in particular useful in the cases close to degenerate`ones when the segments of $E$ tend to vanish or merge.  

\section{Testing the Main formula}
For the verification of the above formula (\ref{main}) we need the good stock of the sets $E$ with analytically known capacities.
We take the inverse polynomial images of segments as such sets. Indeed, let $E:=T^{-1}[-1,1]$, with 
$T=cx^n+lower~terms$ being a degree $n$ polynomial. Then the Green's function has a simple appearance
$$
G_E(x)=\left|
\log|T(x)+\sqrt{T^2(x)-1}|
\right|
$$
wherefrom it immediately follows that $Cap(E)=(2|c|)^{-1/n}$.

Consider Chebyshev polynomials $T_{2n}(x)$ of even degrees $2n=4,6,8$. We compute the capacity of the sets 
$$
E_n:=T^{-1}_{2n}[0,1]=[-1; \cos\frac{\pi}{4n}] \cup \bigcup_{s=1}^{n-1} [\cos\frac{(4s-1)\pi}{4n}; \cos\frac{(4s+1)\pi}{4n}] \cup [\cos\frac{(4n-1)\pi}{4n}; 1].
$$ 
by formula (\ref{main}) and compare it to the analytic value 
\be
\label{CapEn}
Cap(E_n)=2^{-1/n}.
\ee

Below we provide the following data for each $n=2,3,4$: period matrix of the associated curve, image of infinity under the Abel-Jacobi map (up to the apparent accuracy), capacity $C_n$ calculated by formula (\ref{main}) and its comparison with the value (\ref{CapEn}).

{\bf n=2}
$$
Im~\Pi = 
\left(
\begin{array}{ccc}
0.60355339059327 & 0.10355339059327 \\ 0.10355339059327 & 0.60355339059327
\end{array}
\right)
$$
$$
u_{\infty}=\left(\begin{array}{c}
0.37499999999998 \\
0.12499999999992 
\end{array}\right)
$$
$$
C_2=0.84089641525372; \qquad |C_2-Cap(E_2)| < 10^{-14}
$$

{\bf n=3}
$$
Im~ \Pi = 
\left(
\begin{array}{ccc}
0.6220084679281 & 0.1666666666666 & 0.0446581987385 \\
0.1666666666667 & 0.8333333333333 & 0.1666666666667 \\ 
0.0446581987385 &  0.1666666666667 & 0.6220084679281
\end{array}
\right)
$$
$$
u_{\infty}=\left(\begin{array}{ccc}
0.4166666666666 \\
0.2499999999999 \\
0.0833333333333 
\end{array}\right)
$$
$$
C_3= 0.8908987181402;  \qquad |C_3-Cap(E_3)| < 10^{-13}.
$$

{\bf n=4}
$$
Im~ \Pi = 
\left(
\begin{array}{cccc}
0.628417436515 & 0.187075720333 & 0.083522329739 & 0.024864045922 \\ 0.187075720333 & 0.899015486588 & 0.295462095995 & 0.083522329739 \\ 0.083522329739 & 0.295462095995 & 0.899015486588 & 0.187075720333 \\ 0.024864045922 & 0.083522329739 & 0.187075720333 & 0.628417436515
\end{array}
\right)
$$
$$
u_{\infty}=\left(\begin{array}{cccc}
0.437499999999 \\
0.312499999999 \\
0.187499999999 \\
0.062499999999
\end{array}\right)
$$
$$
C_4=0.917004043204 \qquad |C_4-Cap(E_4)| < 10^{-12}
$$

\vspace{5mm}
\parbox{7cm}
{\it
Andrei Bogatyrev\\
Institute for Numerical Mathematics,\\
Russian Academy of Sciences;\\
Moscow Inst. of Physics and Technology;
Moscow State University\\[1mm]
{\tt ab.bogatyrev@gmail.com}}
\hfill
\parbox{7cm}
{\it
Oleg Grigoriev\\
Institute for Numerical Mathematics,\\
Russian Academy of Sciences;\\
{\tt guelpho@mail.ru}}

\end{document}